\magnification\magstep1

\font\ninerm=cmr9    \font\sixrm=cmr6
\font\ninei=cmmi9    \font\sixi=cmmi6
\font\ninesy=cmsy9  \font\sixsy=cmsy6
\font\ninebf=cmbx9  \font\sixbf=cmbx6
\font\nineit=cmti9  
\font\ninett=cmtt9  
\font\ninesl=cmsl9

\font\tengoth=eufm10  \font\ninegoth=eufm9
 \font\sevengoth=eufm7 
\font\sixgoth=eufm6   \font\fivegoth=eufm5
\newfam\gothfam \def\goth{\fam\gothfam\tengoth} 
\textfont\gothfam=\tengoth
\scriptfont\gothfam=\sevengoth 
\scriptscriptfont\gothfam=\fivegoth

\catcode`@=11
\newskip\ttglue

\def\tenpoint{\def\rm{\fam0\tenrm}
  \textfont0=\tenrm \scriptfont0=\sevenrm
  \scriptscriptfont0\fiverm
  \textfont1=\teni \scriptfont1=\seveni
  \scriptscriptfont1\fivei 
  \textfont2=\tensy \scriptfont2=\sevensy
  \scriptscriptfont2\fivesy 
  \textfont3=\tenex \scriptfont3=\tenex
  \scriptscriptfont3\tenex 
  \textfont\itfam=\tenit\def\it{\fam\itfam\tenit}%
  \textfont\slfam=\tensl\def\sl{\fam\slfam\tensl}%
  \textfont\ttfam=\tentt\def\tt{\fam\ttfam\tentt}%
  \textfont\gothfam=\tengoth\scriptfont\gothfam=\sevengoth 
  \scriptscriptfont\gothfam=\fivegoth
  \def\goth{\fam\gothfam\tengoth}
  \textfont\bffam=\tenbf\scriptfont\bffam=\sevenbf
  \scriptscriptfont\bffam=\fivebf
  \def\bf{\fam\bffam\tenbf}%
  \tt\ttglue=.5em plus.25em minus.15em
  \normalbaselineskip=12pt \setbox\strutbox\hbox{\vrule
  height8.5pt depth3.5pt width0pt}%
  \let\big=\tenbig\normalbaselines\rm}

\def\ninepoint{\def\rm{\fam0\ninerm}
  \textfont0=\ninerm \scriptfont0=\sixrm
  \scriptscriptfont0\fiverm
  \textfont1=\ninei \scriptfont1=\sixi
  \scriptscriptfont1\fivei 
  \textfont2=\ninesy \scriptfont2=\sixsy
  \scriptscriptfont2\fivesy 
  \textfont3=\tenex \scriptfont3=\tenex
  \scriptscriptfont3\tenex 
  \textfont\itfam=\nineit\def\it{\fam\itfam\nineit}%
  \textfont\slfam=\ninesl\def\sl{\fam\slfam\ninesl}%
  \textfont\ttfam=\ninett\def\tt{\fam\ttfam\ninett}%
  \textfont\gothfam=\ninegoth\scriptfont\gothfam=\sixgoth 
  \scriptscriptfont\gothfam=\fivegoth
  \def\goth{\fam\gothfam\tengoth}
  \textfont\bffam=\ninebf\scriptfont\bffam=\sixbf
  \scriptscriptfont\bffam=\fivebf
  \def\bf{\fam\bffam\ninebf}%
  \tt\ttglue=.5em plus.25em minus.15em
  \normalbaselineskip=11pt \setbox\strutbox\hbox{\vrule
  height8pt depth3pt width0pt}%
  \let\big=\ninebig\normalbaselines\rm}

\newcount\newpen\newpen=50
\outer\def\beginsection#1\par{\vskip0pt
  plus.3\vsize\penalty\newpen\newpen=-50\vskip0pt
  plus-.3\vsize\bigskip\vskip\parskip
  \message{#1}\leftline{\bf#1}
  \nobreak\smallskip\noindent}
\font\tengoth=eufm10  \font\ninegoth=eufm9
 \font\sevengoth=eufm7 
\font\sixgoth=eufm6   \font\fivegoth=eufm5
\newfam\gothfam \def\goth{\fam\gothfam\tengoth} 
\textfont\gothfam=\tengoth
\scriptfont\gothfam=\sevengoth 
\scriptscriptfont\gothfam=\fivegoth
\def\ZZ{{\bf Z}}
\def\FF{{\bf F}}
\def\QQ{{\bf Q}}
\def\zmod#1{\,\,({\rm mod}\,\,#1)}

\def\bibliography#1\par{\vskip0pt
  plus.3\vsize\penalty-200\vskip0pt
  plus-.3\vsize\bigskip\vskip\parskip
  \message{Bibliography}\leftline{\bf
  Bibliography}\nobreak\smallskip\noindent
  \ninepoint\frenchspacing#1}
\outer\def\beginsection#1\par{\vskip0pt
  plus.3\vsize\penalty\newpen\newpen=-50\vskip0pt
  plus-.3\vsize\bigskip\vskip\parskip
  \message{#1}\leftline{\bf#1}
  \nobreak\smallskip\noindent}
\def\today{\ifcase\month\or January \or February\or
March\or April\or May\or June\or July\or August\or
September\or October\or November\or December\fi
\space\number\day, \number\year}

\def\APR{{1}}
\def\AH{{2}}
\def\AKS{{3}}
\def\AR{{4}} 
\def\AM{{5}}
\def\BA{{6}} 
\def\COH{{7}}
\def\CP{{8}} 
\def\GK{{9}} 
\def\KP{{10}}
\def\LANG{{11}}
\def\LENMIL{{12}}
\def\LEN{{13}}
\def\LRES{{14}}
\def\ECM{{15}}
\def\LP{{16}}
\def\MIH{{17}}
\def\MIL{{18}}
\def\MOR{{19}}
\def\RAB{{20}}
\def\SCH{{21}} 
\def\WASH{{22}}

\noindent Ren\'e Schoof\hfill  Amsterdam, \today
\smallskip
\hrule
\bigskip\bigskip\bigskip\bigskip
\centerline{\bf Four primality testing algorithms}
\bigskip\bigskip
\beginsection Introduction.

In this expository paper we describe four primality tests. The first test is very efficient, but is only capable of proving  that
a given number is either composite or `very probably' prime. The second test is a deterministic polynomial time algorithm to prove that a given numer is either prime or composite. The third and fourth primality tests are at present  most widely used in practice. Both tests are capable of proving that a given number is prime or composite, but neither algorithm is deterministic.
The third algorithm exploits the arithmetic of cyclotomic fields. Its running time is almost, but not quite polynomial time. The
fourth algorithm exploits elliptic curves. Its running time is difficult to estimate, but  it behaves well in practice.

In section~1 we discuss the Miller-Rabin test. This is one of the most efficient probabilistic primality tests.
Strictly speaking, the Miller-Rabin test is not a primality test but rather a `compositeness test', since it does not prove the primality of a number.   Instead, if $n$ is {\it not} prime, the algorithm proves this in all likelihood very quickly.
On the other hand, if $n$ happens to be prime,
the algorithm merely provides strong evidence for its primality.  Under the assumption of the Generalized Riemann Hypothesis one can turn the Miller-Rabin algorithm into a deterministic polynomial time primality test. This idea, due to G.~Miller, is also explained.

In section 2 we describe the deterministic polynomial time  primality test~[\AKS]  that was proposed by M.~Agrawal, N.~Kayal and N.~Saxena in~2002.  At the moment the present paper was written, this new test, or rather a more efficient probabilistic version of it, had not yet been widely implemented.
In practice, therefore, for {\it proving} the primality of a given integer, 
one still relies on older tests that are either not provably polynomial time or not deterministic.
In the remaining two sections we present the two most widely used such tests. 

In section 3 we discuss the cyclotomic primality test. 
This test is deterministic and  is actually capable of {\it proving} that a given integer $n$ is either prime or composite. It  does not run in polynomial time, but very nearly so. We describe a practical non-deterministic version of the algorithm.
Finally in section~4, we describe the elliptic curve primality test. This algorithm 
also provides a {\it proof} of the primality or compositeness of a given integer~$n$. Its running time is hard to analyze, but in practice the algorithm  seems to run in polynomial time. It is not deterministic. The
 two `practical' tests described in sections~3 and~4 have been implemented and fine tuned. Using either of them it is now possible to routinely prove the primality of numbers that have  several thousands of decimal digits~[\MIH, \MOR].

\beginsection 1. A probabilistic  test.

In this section we present a practical and efficient probabilistic primality test.
Given a composite integer $n>1$, this algorithm proves with high probability very 
quickly that $n$ is not prime. On the other hand, if $n$ passes the
test, it is merely {\it likely} to be prime. The algorithm consists of repeating one simple 
step, a Miller-Rabin test, several times with different random initializations.
The probability that  a composite number
is {\it not} recognized as such by the algorithm, can be made arbitrarily small
 by repeating the main step  a number of times. 
 The algorithm was first proposed by~M.~Artjuhov~[\AR] in~1966. 
 In 1976 M.~Rabin proposed the probabilistic version~[\RAB].
Under assumption of the Generalized Riemann Hypothesis (GRH) one can actually {\it prove} that $n$ is prime by 
applying the test sufficiently often.
This leads to G.~Miller's  {\it conditional}  algorithm~[\MIL]. Under assumption of~GRH it runs in polynomial time.
Our presentation  follows the  presentation of the algorithms in the excellent book by
R.~Crandall and C.~Pomerance~[\CP].

The following theorem is the key ingredient.

\proclaim Theorem 1.1. Let $n>9$ be an odd positive {\it composite} integer. We write $n-1=2^km$ for some
exponent $k\ge 1$ and some odd integer~$m$.
Let
$$
B=\{x\in(\ZZ/n\ZZ)^*: \hbox{$x^m=1$ or $x^{m2^i}=-1$ for some $0\le i<k$}\}.
$$
Then we have
$$
{{\#B}\over{\varphi(n)}}\,\,\le\,\, {1\over 4}.
$$
Here  $\varphi(n)=\#(\ZZ/n\ZZ)^*$ denotes Euler's $\varphi$-function.

\smallskip\noindent{\bf Proof.} Let $2^l$ denote the largest power of 2 that has the property
that it divides $p-1$ for 
every prime $p$ divisor of~$n$. Then the set $B$ is contained in
$$
B'=\{x\in(\ZZ/n\ZZ)^*: \,x^{m2^{l-1}}=\pm1\}.
$$
Indeed, clearly any $x\in(\ZZ/n\ZZ)^*$ satisfying $x^m=1$ is contained in~$B'$. On the other hand,
if $x^{m2^i}=-1$ for some $0\le i<k$, we have $x^{m2^i}\equiv-1\zmod{p}$
for every prime $p$ dividing~$n$.
It follows that for every~$p$, the {\it exact} power of 2 dividing 
the order of~$x$ modulo~$p$,   is equal to~$2^{i+1}$.
In particular, $2^{i+1}$ divides $p-1$ for  every prime divisor $p$ of~$n$. Therefore we have $l\ge i+1$.
So we can write  that $x^{m2^{l-1}}=(-1)^{2^{l-i-1}}$, which  is~$-1$ or~$+1$ depending on whether $l=i+1$ or~$l>i+1$. 
It follows that $B\subset B'$.

By the Chinese Remainder Theorem, the number of elements $x\in(\ZZ/n\ZZ)^*$ for which we have  $x^{m2^{l-1}}=1$,
is equal to the product over $p$ of the number of solutions to the equation $X^{m2^{l-1}}=1$ modulo~$p^{a_p}$.
Here  $p$ runs over the prime divisors of~$n$ and $p^{a_p}$ is the exact power of~$p$ dividing~$n$.
Since each of the groups $(\ZZ/p^{a_p}\ZZ)^*$ is cyclic, the number of solutions modulo $p^{a_p}$ is given by ${\rm gcd}((p-1)p^{a_p-1},m2^{l-1})={\rm gcd}(p-1,m)2^{l-1}$. The last equality follows from the fact that $p$ does not divide~$m$.
Therefore we have
$$
\#\{x\in(\ZZ/n\ZZ)^*: \,x^{m2^{l-1}}=1\}\quad=\quad\prod_{p|n}{\rm gcd}(p-1,m)2^{l-1}.
$$
Similarly, the number of solutions of the equation $X^{m2^{l}}=1$ modulo~$p^{a_p}$ is
equal to ${\rm gcd}(p-1,m)2^{l}$, which is twice the number of solutions 
of~$X^{m2^{l-1}}=1$ modulo~$p^{a_p}$. It follows that the number of solutions of
the equation $X^{m2^{l-1}}=-1$ modulo~$p^{a_p}$ is also equal to ${\rm gcd}(p-1,m)2^{l-1}$.
Therefore we have
$$
\#B'\quad=\quad2\prod_{p|n}{\rm gcd}(p-1,m)2^{l-1},
$$
and hence
$$
{{\#B'}\over{\varphi(n)}}\quad=\quad2\prod_{p|n}{{{\rm gcd}(p-1,m)2^{l-1}}\over{(p-1)p^{a_p-1}}}.
$$
Suppose now that the propertion ${{\#B}\over{\varphi(n)}}$ {\it exceeds}  ${1\over 4}$. 
We want to derive a contradiction.
Since we have $B\subset B'$, the
inequality above  implies that
$$
{1\over 4}\quad<\quad 2\prod_{p|n}{{{\rm gcd}(p-1,m)2^{l-1}}\over{(p-1)p^{a_p-1}}}.\eqno{(*)}
$$
We draw a number of conclusions from this inequality. First we 
note that ${\rm gcd}(p-1,m)2^{l-1}$ divides $(p-1)/2$ so that the right hand side of $(*)$ is at most $2^{1-t}$ where $t$ is the number of different primes dividing~$n$. It follows that $t\le 2$. 

Suppose that  $t=2$, so that $n$ has precisely two distinct prime divisors.
If one of them, say $p$, has the property that $p^2$ divides $n$ so that $a_p\ge 2$, then the right hand side of $(*)$ is at most $2^{1-2}/3=1/6$. Contradiction. It follows that all exponents $a_p$ are equal to~1, so that
$n=pq$ for two distinct primes $p$ and~$q$. The inequality ($*$) now becomes 
$${{p-1}\over{{\rm gcd}(p-1,m)2^{l}}}\cdot {{q-1}\over {{\rm gcd}(q-1,m)2^{l}}}<2.
$$
Since the factors on the left hand side of this inequality are positive integers, they are both equal to~1. This implies that
$p-1={\rm gcd}(p-1,m)2^{l}$ and $q-1={\rm gcd}(q-1,m)2^{l}$. It follows that the exact power of 2 dividing  $p-1$ as well as the exact power of 2 dividing $q-1$ are equal to~$2^l$ and that the odd parts of $p-1$ and $q-1$ divide~$m$. 
Considering the relation $pq=1+2^km$ modulo the odd part of~$p-1$, we see that the odd part of $p-1$ divides the odd part of~$q-1$. By symmetry, the odd parts of $p-1$ and $q-1$ are therefore equal.
This implies $p-1=q-1$ and contradicts the fact that $p\not=q$. Therefore we have $t=1$ and hence $n=p^a$ for some odd prime~$p$ and exponent~$a\ge 2$. The inequality ($*$) now says that $p^{a-1}<4$, so that $p=3$ and $a=2$, contradicting the hypothesis that~$n>9$. This proves the Theorem.

\bigskip
When a random  $x\in(\ZZ/n\ZZ)^*$ is checked to be contained in the set $B$ of   Theorem~1.1, 
we say that `$n$ passes a Miller-Rabin test'.
Checking that $x\in B$ involves raising  $x\in\ZZ/n\ZZ$ to an exponent 
that is no more than~$n$. Using the binary expansion of the exponent, this takes no more that $O(\log\,n)$ multiplications
in~$\ZZ/n\ZZ$. Therefore a single exponentiation involves $O((\log\,n)^{1+\mu})$ elementary operations or bit operations. Here $\mu$ is a constant with the property that the multiplication algorithm in~$\ZZ/n\ZZ$ takes
no more than $O((\log\,n)^{\mu})$ elementary operations. We have that $\mu=2$ when we use the usual multiplication algorithm, while one can take $\mu=1+\varepsilon$ for any $\varepsilon>0$ by employing fast multiplication techniques.

By Theorem~1.1 the probability that a composite number $n$ passes a single Miller-Rabin test, is at most 25\%.
Therefore, the probability that $n$ passes $\log\,n$ such tests is smaller than $1/n$. 
The probability that a large composite $n$ passes $(\log\,n)^2$ tests is astronomically small: less than $n^{-\log\,n}$.
Since for most composite~$n$ the probability that $n$ passes a Miller-Rabin test is much smaller than $1/4$, 
one is in practice already convinced of the primality of~$n$,
when $n$ successfully passes a handful of Miller-Rabin tests. This is enough for most commercial applications.

\hyphenation{Di-ri-chlet}
Under assumption of the Generalized Riemann Hypothesis (GRH) for quadratic Dirichlet characters, the Miller-Rabin test can be transformed into a  {\it deterministic} polynomial time primality test. This result goes back to G.~Miller~[\MIL].

\proclaim Theorem 1.2. (GRH) Let $n$ be an odd positive {\it composite} integer. Let  $n-1=2^km$ for some exponent $k\ge 1$ and some odd integer~$m$. If for all integers $x$ between 1 and $2(\log\,n)^2$ one has
$$
\hbox{$x^m\equiv1\zmod n$ \quad or \quad  $x^{2^im}\equiv-1\zmod n$ for some $0\le i<k$},
$$
then $n$ is a prime number.

\smallskip\noindent{\bf Proof.} We first show that $n$ is squarefree. See also~[\LENMIL].
Suppose that $p$ is a prime for which $p^2$ divides~$n$. A special case of a result of Konyagin and Pomerance~[\KP,~(1.45)] on the distribution of smooth numbers implies that  for every odd integer $r\ge 5$ 
one has that
$$
\#\{a\in\ZZ:\hbox{$1\le a\le r$ and $a$ is product of primes $\le(\log\,r)^2$}\}\ \ \ge\ \ \sqrt{r}.
$$
We apply this with $r=p^2$. It follows that the subgroup $H$ of $(\ZZ/p^2\ZZ)^*$ that is generated by
the natural numbers $x\le (\log\,n)^2$ has order at least~$p$. On the other hand, the hypothesis of the theorem
implies that every $x\in H$, being   a product of numbers $a$ that satisfy  $a^{n-1}\equiv 1\zmod{p^2}$, 
 satisfies  $x^{n-1}
\equiv 1\zmod{p^2}$. Since the order of the group $(\ZZ/p^2\ZZ)^*$ is $p(p-1)$ and $p$ does not divide~$n-1$,
we see that any $x\in H$ must satisfy  $x^{p-1}\equiv 1\zmod{p^2}$. But this is impossible, because the subgroup 
 of $(\ZZ/p^2\ZZ)^*$ that consists of  elements having this property, has order~$p-1$. 

Therefore, if $n$ is composite, it is  divisible by two  odd distinct primes $p$ and~$q$. Let $\chi$  denote the quadratic character of conductor $p$. By a result of E.~Bach~[\BA], {\it proven under assumption of the GRH}, there exists a natural number $x\le 2(\log\,p)^2 < 2(\log\,n)^2$ for which~$\chi(x)\not=1$. Since the condition of the theorem implies that we have ${\rm gcd}(x,n)=1$, we must have~$\chi(x)=-1$.
Writing $p-1=2^l\mu$ for some exponent $l\ge1 $ and some 
odd integer~$\mu$, we have that $x^{2^{l-1}\mu}\equiv\chi(x)=-1\zmod{p}$. 
This implies that $-1$ is contained in the subgroup of $(\ZZ/p\ZZ)^*$ generated by~$x$. 
Since the 2-parts of the subgroups of $(\ZZ/p\ZZ)^*$ generated by $x^m$ and by $x$  are the same, we have $x^m\not\equiv 1\zmod p$ and hence~$x^m\not\equiv1\zmod{n}$.
Therefore  the hypothesis of the theorem implies that $x^{2^im}\equiv-1\zmod n$ for some $0\le i<k$.  
Since for this value of $i$ we also have $x^{2^im}\equiv-1\zmod p$,  necessarily the equality $i=l-1$ holds. 
It follows that we have $x^{2^{l-1}m}\equiv-1\zmod q$, so that the order of $x^m\zmod q$ is equal to~$2^l$. Writing 
$q-1=2^{l'}\mu'$ for some exponent $l'\ge 1$ and some odd integer~$\mu'$, we have therefore~$l\le l'$.

Repeating the argument, but switching the roles of $p$ and~$q$, we conclude that $l=l'$. Let $\chi'$ denote  the quadratic character of conductor $q$. A second application of Bach's theorem, this time to the {\it non-trivial} character~$\chi\chi'$, provides us with  a natural number $y\le 2(\log\,n)^2$ for which $\chi\chi'(y)\not=1$ and hence, say, $\chi(y)=-1$ while $\chi'(y)=1$. The arguments given above, but this time applied to~$y$, show that 
we cannot have $y^{m}\equiv-1\zmod n$, so that  necessarily 
$y^{2^im}\equiv-1\zmod n$ for some $0\le i<k$. Moreover, the exponent $i$ is equal to $l-1=l'-1$. 
It follows that $y^{2^{l'-1}m}\equiv-1\zmod q$. This implies that the element $y^m\in(\ZZ/q\ZZ)^*$ has order $2^{l'}$.
Since the subgroups of $(\ZZ/q\ZZ)^*$ generated by $y^m$ and $y^{\mu'}$  are equal, the order 
of $y^{\mu'}\in(\ZZ/q\ZZ)^*$ is also $2^{l'}$. This contradicts the fact that 
$1=\chi'(y)\equiv y^{2^{l'-1}\mu'}\zmod{q}$.

We conclude that $n$ is prime and the result follows.
\bigskip
It is clear how to apply Theorem~1.2  and obtain a  test that proves that~$n$ is prime under condition of~GRH: given an odd integer $n>1$, we simply test the condition of Theorem~1.2 for all $a\in\ZZ$ satisfying~$1<a<2(\log\,n)^2$.  If $n$ passes all these tests and GRH holds, then $n$ is prime.  Each test involves an exponentiation in the ring $\ZZ/n\ZZ$. Since the exponent is less than~$n$, this can be done using only $O((\log\,n)^{1+\mu})$ elementary operations. Therefore this is a polynomial time primality test. Testing~$n$ takes $O((\log\,n)^{3+\mu})$ elementary operations. As before,
we have  $\mu=2$ when we use the usual multiplication algorithm, while we can take $\mu=1+\varepsilon$ for any $\varepsilon>0$ by employing fast multiplication techniques.

\beginsection 2. A deterministic polynomial time primality test.

In the summer of~2002 the three Indian computer scientists 
M.~Agrawal, N.~Kayal and N.~Saxena presented a deterministic polynomial time primality test. We describe and analyze 
this extraordinary result in this section.

 For any prime number $r$ we let
$\Phi_r(X)=X^{r-1}+\ldots+X+1$ denote the $r$-th cyclotomic polynomial. Let~$\zeta_r$ be a zero of~$\Phi_r(X)$
and let $\ZZ[\zeta_r]$ denote the ring generated by~$\zeta_r$ over~$\ZZ$. For any $n\in\ZZ$ we write
$\ZZ[\zeta_r]/(n)$ for the residue ring $\ZZ[\zeta_r]$ modulo the ideal $(n)$ 
generated by~$n$. For $n\not=0$, this is a finite ring.

\proclaim Theorem 2.1.  Let $n$ be an odd positive integer and let $r$ be a prime number. Suppose that
\smallskip\item{(i)} $n$ is not divisible by any of the primes $\le r$;
\smallskip\item{(ii)} the order of $n\zmod r$ is at least $(\log\,n/\log\,2)^2$;
\smallskip\item{(iii)} for every $0\le j< r$ we have $(\zeta_r+j)^n=\zeta_r^n+j$ in~$\ZZ[\zeta_r]/(n)$.
\smallskip\noindent
Then $n$ is a prime power.

\smallskip\noindent{\bf Proof.}  It follows from condition~{\sl (ii)} that we have $n\not\equiv 1\zmod r$. Therefore there exists a prime divisor $p$ of~$n$ that is not congruent to~$1\zmod r$. Let $A$ denote the  $\FF_p$-algebra~$\ZZ[\zeta_r]/(p)$. It is a quotient of the ring~$\ZZ[\zeta_r]/(n)$. For 
$k\in\ZZ$ coprime to~$r$ we let $\sigma_k$ denote the ring automorphism of~$A$ determined by $\sigma_k(\zeta^{}_r)=\zeta_r^k$. The map $(\ZZ/r\ZZ)^*\mapsto \Delta$ given by $k\mapsto \sigma_k$ is a well defined isomorphism. We single out two special elements of~$\Delta$. One is the {\it Frobenius automorphism}  $\sigma_p$ and the other is~$\sigma_n$. Let $\Gamma$ denote the subgroup of~$\Delta$ that is generated by~$\sigma_p$ and~$\sigma_n$. 

Next we consider the subgroup $G$ of elements of the multiplicative group 
$A^*$ that are annihilated by the endomorphism $\sigma_n-n\in\ZZ[\Delta]$. In other words, we put
$$
G=\{a\in A^*: \sigma_n(a)=a^n\}.
$$
Pick a maximal ideal ${\goth m}$ of~$A$  and put $k=A/{\goth m}$. Then $k$ is a finite extension of~$\FF_p$,  generated by a primitive $r$-th root of unity.
Let $H\subset k^*$ be the image of $G$ under the natural map~$\pi:A\longrightarrow k$. 
The group $H$ is cyclic. Its order is denoted by~$s$.
We have the following commutative diagram.
$$
\def\normalbaselines{\baselineskip20pt\lineskip3pt 
\lineskiplimit3pt}
\matrix{G&\subset &A^*\cr
\downarrow\rlap{\hbox{$\scriptstyle
\pi$}}&&\downarrow\rlap{\hbox{$\scriptstyle \pi$}}\cr
H&\subset& k^*&\cr}
$$
Since $\Delta$ is commutative, it acts on~$G$. Since $\sigma_n$ and 
$\sigma_p$ act on $G$ by raising to the power $n$ and $p$ respectively, 
every $\sigma_m\in\Gamma$  acts by raising $g\in G$ to a certain power $e_m$
that is prime to~$\#G$. 
The powers $e_m$ are well determined modulo the exponent ${\rm exp}(G)$ of $G$.
Therefore the map $\Gamma\longrightarrow(\ZZ/{\rm exp}(G)\ZZ)^*$,
given by $\sigma_m\mapsto e_m$, is a well defined group homomorphism.
Since $H$ is a cyclic quotient of~$G$, its order $s$ divides the exponent of~$G$ and the map  
$\sigma_m\mapsto e_m$ induces a homomorphism
$$
\Gamma\quad\longrightarrow\quad(\ZZ/s\ZZ)^*.
$$
If $m\equiv p^in^j\zmod r$, then it maps $\sigma_m\in\Gamma$ to $e_m\equiv p^in^j\zmod{s}$.

It is instructive to see what  all this boils down to when $n$ is prime. 
Then  we have $n=p$ and~$\sigma_n$ is equal to the Frobenius automorphism~$\sigma_p$. The group
$G$ is all of~$A^*$ so that~$H$ is equal to~$k^*$. Writing $f$ for the order of $p$ modulo~$r$, the
group $\Gamma=\langle\sigma_p\rangle$ has order~$f$ while
the groups $H=k^*$ and its automorphism group ${\rm Aut}(H)$ are {\it much larger}.
Indeed,  $H$ has order
$s=p^f-1=n^{\#\Gamma}-1$ and ${\rm Aut}(H)\cong(\ZZ/s\ZZ)^*$ is of comparable size

Under the conditions of the theorem, but {\it without assuming} that $n$ is prime, something
similar can be shown to be true.
 
\smallskip\noindent {\bf Claim.} We have that
$$
s\,\,>\,\,n^{\left[\sqrt{\#\Gamma}\right]}_{}.
$$
Using this inequality, we complete the proof of the theorem.
Consider the homomorphism
$$\Gamma\,\,\longrightarrow\,\,(\ZZ/s\ZZ)^*
$$
constructed above. We first  apply the box principle in the
{\it small} group $\Gamma$  and then obtain a relation in $\ZZ$ from a relation in
$(\ZZ/s\ZZ)^*$ using the fact that the latter group is {\it very large}. 

Let  $q=n/p$. We consider the products
$\sigma_p^i\sigma_q^j\in\Gamma$ for $0\le i,j\le \left[\sqrt{\#\Gamma}\right]$. Since we have
$(1+\left[\sqrt{\#\Gamma}\right])^2>\#\Gamma$, there are two pairs $(i,j)\not=(i',j')$
for which  $\sigma_p^i\sigma_q^j$ and $\sigma_p^{i'}\sigma_q^{j'}$ are the same element
in~$\Gamma$.  It follows that their images in the
group $(\ZZ/s\ZZ)^*$ are the same as well. 
Since $\sigma_q$ is mapped to~$q\zmod s$,  this means that $p^i{q}^j\equiv p^{i'}q^{j'}\zmod{s}$.
The integer
$p^i{q}^j$ does not exceed
$n^{{\rm max}(i,j)}\le n^{\left[\sqrt{\#\Gamma}\right]}< s$.
The same holds for $p^{i'}q^{j'}$. We conclude that $p^i{q}^j=p^{i'}q^{j'}$ in~$\ZZ$!
Since $(i,j)\not=(i',j')$ it follows that $n$ is a power of~$p$.

This proves the theorem.
\bigskip
\noindent{\bf Proof of the claim.}  We first estimate $s=\#H$  in terms of $\#G$. Then we show that $G$ is large. 

The first bound we show is
$$
s\ge \#G^{1/[\Delta:\Gamma]}.\eqno{(*)}
$$
Let $C$ denote a set of coset representatives of $\Gamma$ in~$\Delta$ and consider the  homomorphism
$$
G\longrightarrow\prod_{i\in C}k^*
$$
given by mapping $a\in G$ to the vector $(\sigma_i(a)\zmod{{\goth m}})_{i\in C}$. 

This map is injective. Indeed, if $a\in G$
has the property that $\sigma_i(a)=1$ for some~$i$, then we also have 
$\sigma_{in}(a)=\sigma_i(a^n)=\sigma_i(a)^n=1$ and similarly~$\sigma_{ip}(a)=1$.
In other words, we have $\sigma(a)=1$ for all elements $\sigma$ in the coset of~$\Gamma$ containing~$\sigma_i$.
Therefore, if $a\in G$ has the property that $\sigma_i(a)=1$ for all~$i\in C$, then automatically also
$\sigma_i(a)=1$ for all~$i\in(\ZZ/r\ZZ)^*$. It follows that $\sigma_i(a-1)=0$ for all~$i\in(\ZZ/r\ZZ)^*$. Writing the element $a-1$ as $f(\zeta_r)$ for some polynomial~$f(X)\in\FF_p[X]$, this implies that $f(\zeta_r^i)=0$ for all~$i\in(\ZZ/r\ZZ)^*$. It follows that the cyclotomic polynomial $\Phi_r(X)$ divides~$f(X)$ in $\FF_p[X]$ and hence that $a-1=0$, as required.

Since for every $i\in C$, the image of the map $G\longrightarrow k^*$ given by 
$a\mapsto \sigma_i(a)\zmod{\goth m}$ is equal to~$H$, 
the  injectivity of the homomorphism implies that
$\#G\le s^{[\Delta:\Gamma]}$ as required.

\smallskip
The second estimate is
$$
\#G\ge 2^{r-1}.\eqno{(**)}
$$
Since we have $p\not\equiv1\zmod r$, the irreducible factors of 
$\Phi_r(X)=(X^r-1)/(X-1)$ in the ring $\FF_p[X]$  have degree at
least~2 and hence cannot divide any polynomial of degree~1. Therefore the elements
$\zeta_r+j$   for $0\le j< r-1$ are not contained in any maximal
ideal of the ring~$A$. It follows that  they are {\it units} of~$A$. 
By condition~{\sl (iii)},  for each subset $J\subset\{0,1,\ldots,r-2\}$ the element  
$$\prod_{j\in J}(\zeta_r^{}+j)
$$ 
is contained in~$G$. 

All these elements are {\it distinct}.  Indeed, since the degree of the cyclotomic polynomial $\Phi_r$ is $r-1$, the only two elements that could be equal are the ones corresponding to the extreme cases
$J=\emptyset$ and to $J=\{0,1,\ldots,r-2\}$. This can only happen when
$\prod_{j=0}^{r-2}(X+j)-1$ is divisible by $\Phi_r(X)$ in the ring $\FF_p[X]$. 
Since both polynomials have the same degree,
we then necessarily  have
$\prod_{j=0}^{r-2}(X+j)-1=\Phi_r(X)$.
Inspection of the constant terms shows that $p=2$. But this is impossible, because  $n$ is odd. 

Since there are $2^{r-1}$ subsets  $J\subset\{0,1,\ldots,r-2\}$, we conclude that 
$\#G\ge 2^{r-1}$.
as required.

\smallskip
Combining the  inequalities ($*$) and ($**$) we find that
$$
s\,\,\ge\,\,
\#G_{}^{{1/{[\Delta:\Gamma]}}}\,\,\ge\,\,2_{}^{(r-1)/
{[\Delta:\Gamma]}}\,\,=\,\,2^{\#\Gamma}
\,\,>\,\,
 n^{\sqrt{\#\Gamma}}\,\,\ge\,\, n^{\left[\sqrt{\#\Gamma}\right]}.
$$
Here we used the inequality $\#\Gamma>(\log
n/\log 2)^2$. It follows from the fact that the order of $\sigma_n\in\Gamma$
 is larger than~$(\log n/\log 2)^2$. Indeed, this order is equal to
 the order of $n$ modulo~$r$, which  by condition~{\sl (ii)}
 is larger than~$(\log n/\log 2)^2$. 

This proves the claim.
\medskip
\bigskip\noindent This theorem leads to the following primality 
test.  

\bigskip

\proclaim Algorithm 2.2. Let $n> 1$ be a given odd integer. 
\smallskip
{\item{(i)} First check that $n$ is not a proper power of an integer.
\smallskip
\item{(ii)} By successively trying $r=2,3,\ldots$, determine the
smallest prime $r$ not dividing $n$ nor any of the numbers $n^i-1$ for
$0\le i\le (\log n/\log 2)^2$.\smallskip
\item{(iii)} For $0\le j< r-1$ check that $ (\zeta_r+j)^n= \zeta_r^n +j$ in the
ring $\ZZ[\zeta_r]/(n)$.}
\smallskip\noindent
If the number $n$ does not pass the tests, it is composite.
If it passes them, it is a prime.

\smallskip\noindent{\bf Proof of correctness.} If $n$ is prime,
it  passes the tests  by Fermat's little theorem.
Conversely, suppose that $n$ passes the tests. 
We check the conditions
of Theorem~2.1.  By definition of~$r$, the number~$n$ has no prime divisors $\le r$. Since
$r$ does not divide any of the $n^i-1$ for $1\le i\le(\log n/\log 2)^2$, the
order of $n$ modulo $r$ exceeds $(\log n/\log 2)^2$.
This shows that the second condition of Theorem~2.1 is satisfied. Since test~{\sl (iii)} has been 
passed successfully, the third condition is satisfied. We deduce that $n$ is a prime power.
Since $n$ passed the first test, it is therefore prime.

\medskip\noindent{\bf Running time analysis.}
The first test is performed by  checking that $n^{1/m}\not\in\ZZ$  for all integers
$m$ between 2 and~$\log n/\log 2$. This can be done in time~$O((\log n)^4)$
by computing sufficiently  accurate approximations to $n^{1/m}\in{\bf R}$.  
The second test does not take more than $r$ times $O((\log n)^2)$
multiplications with modulus $\le r$.
This takes at most $O(r(\log r\log n)^2)$ bit operations.
The third test takes $r$ times $O(\log n)$ multiplications in the ring 
$\ZZ[\zeta_r]/(n)$.  The latter ring is isomorphic to $\ZZ[X]/(\Phi_r(X),n)$.
If the multiplication algorithm that we use to multiply
two elements of bit size $t$ takes no more than $O(t^{\mu})$ elementary
operations, then this adds up to
$O((r\log n)^{1+\mu})$ elementary operations. 
Since $\mu\ge 1$ and since  $r$ exceeds the order of $n$ modulo $r$, we have $r>(\log n/\log 2)^2$.
Therefore the third test is the dominating part of the algorithm.

We estimate how small we
can take~$r$. By definition of~$r$, the product $n\prod_{i}(n^i-1)$ is  divisible
by all primes $l<r$. Here the product runs over $i\le (\log n/\log 2)^2$.
So
$$
\sum_{l<r}\log l\,\,\le\,\, \log n+\log n\sum_{1\le i\le ({{\log n}\over{\log 2}})^2}i\,\,=\,\, O((\log n)^5).
$$
A weak and easily provable form of the prime number theorem says that
there exists a constant $c>0$, so that for every $r$ we have~$\sum_{l<r}\log l\ge c r$. 
Therefore we have~$r=O((\log n)^5)$.
It follows that the algorithm takes
$O((\log n)^{6(1+\mu)})$ elementary operations.  When the usual
multiplication algorithm is used, we have that $\mu=2$ and this leads to an
algorithm that takes at most
$O((\log n)^{18})$ elementary operations . It takes 
$O((\log n)^{12+\varepsilon})$  elementary operations   when fast
multiplication techniques are employed.

\bigskip\noindent{\bf Remark 1.} Since the upper bound
$\sqrt{\#\Gamma}$ is optimal for the box principle, the inequality
$2^{\#\Gamma}\,>\, n^{\sqrt{\#\Gamma}}$ used above
implies that  $\#\Gamma=r-1$  needs to be at least~$(\log n/\log 2)^2$. 
This we know to be the case   because the order of $\sigma_n\in\Gamma$, which 
is equal to the order of  $n\in(\ZZ/r\ZZ)^*$, exceeds~$(\log n/\log 2)^2$. 
The argument involving the
prime number theorem given above implies then that we cannot expect to
be able to prove that the order of magnitude of the prime  $r$
is smaller than~$O((\log n)^5)$. Therefore 
this algorithm cannot be expected to be proved to run faster than $O((\log n)^{6(1+\mu)})$. 
On the other hand, in practice one easily finds a suitable prime 
of the smallest possible size $O((\log n)^2$. Therefore the
practical running time of the algorithm is~$O((\log n)^{3(1+\mu)})$.

\medskip\noindent{\bf Remark 2.}  One may replace the
ring $\ZZ[\zeta_r]/(n)\cong (\ZZ/n\ZZ)[X]/(\Phi_r(X))$ by
any  Galois extension of $\ZZ/n\ZZ$ of the form $(\ZZ/n\ZZ)[X]/(f(X))$ that admits
an automorphism $\sigma$ with the properties that \smallskip
\item{--} $\sigma(X)=X^n$; \smallskip
\item{--} $\sigma$ has order at least~$(\log n/\log
2)^2$. \smallskip

\noindent
This was pointed out by Hendrik Lenstra shortly after the algorithm described above
came out.
The running time of the resulting modified algorithm is then $O((d\log n)^{1+\mu})$ where $d$ is
the degree of the polynomial~$f(X)$.  Since the order of $\sigma$ is at most~$d$, one has that $d>(\log
n/\log 2)^2$ and one cannot obtain an algorithm that runs faster than $O((\log
n)^{3(1+\mu)})$.  Since then Lenstra and Pomerance~[\LP] showed that for every $\varepsilon>0$
one can construct suitable rings with 
$d=O((\log n)^{2+\varepsilon})$. This leads to a primality test that runs in 
time~$O((\log n)^{(3+\varepsilon)(1+\mu)})$. This is essentially the same as  the practical running time
mentioned above.

\beginsection 3. The cyclotomic primality test. 

In this section  we describe the cyclotomic primality test. This algorithm was proposed in~1981 
by L.~Adleman, C.~Pomerance and R.~Rumely~[\APR]. It is one of the most powerful practical tests available today~[\MIH].
Our exposition  follows H.~Lenstra's Bourbaki lecture~[\LEN]. See also~[\COH, section 9.1] and~[\WASH, section 16.1].
The actual computations involve {\it Jacobi sums}, but the basic idea of the algorithm is best explained in terms of {\it Gaussian sums}. See the books by  L.~Washington~[\WASH] and S.~Lang~[\LANG] for a more systematic discussion of the basic properties of Gaussian sums and Jacobi sums. For any positive integer $r$, we denote the subgroup of $r$-th roots of unity
of $\overline{\QQ}^*$ by~$\mu_r^{}$.
\medskip

\noindent{\bf Definition.} {\sl Let $q$ be a prime and let $r$ be a positive integer prime to~$q$.
Let $\chi:(\ZZ/q\ZZ)^*\longrightarrow \mu_{r}$ be
a  character and let $\zeta_q$ be a primitive $q$-th root of unity.
Then we define the {\it Gaussian sum} $\tau(\chi)$ by
$$
\tau(\chi)=-\sum_{x\in(\ZZ/q\ZZ)^*}\zeta^x_q\chi(x).
$$}
The Gaussian sum $\tau(\chi)$ is an algebraic integer, contained in the  cyclotomic field 
$\QQ(\zeta_r^{},\zeta_q^{})$. 
We have the following diagram of fields
\medskip
$$
\matrix{&&\QQ(\zeta_r,\zeta_q)&&\cr
&&&&\cr
&\hbox{\llap{$G$\ \ }$\nearrow$}&&\hbox{$\nwarrow$\rlap{\ \ $\Delta$}}&\cr
&&&&\cr
\QQ(\zeta_q^{})&&&&\QQ(\zeta_r^{})\cr
&&&&\cr
&\nwarrow&&\nearrow&\cr
&&&&\cr
&&\QQ&&\cr}
$$\medskip\noindent
The Galois group  of $\QQ(\zeta_r^{},\zeta_q^{})$ over $\QQ$  is isomorphic to~$\Delta\times G$. Here we have $\Delta=\{\sigma_i:i\in(\ZZ/r\ZZ)^*\}$,
where  $\sigma_i\in\Delta$ is the automorphism that acts trivially on $q$-th roots of unity, while its action of
$r$-th roots of unity is given by  $\sigma_i(\zeta_{r})=\zeta_{r}^i$. The map $(\ZZ/r\ZZ)^*\longrightarrow\Delta$ given by $i\mapsto \sigma_i$ is an isomorphism of groups.
Similarly, we have 
 $G=\{\rho_j:j\in(\ZZ/q\ZZ)^*\}$
 where  $\rho_j\in\Delta$ is the automorphism given by $\rho_j(\zeta_r)=\zeta_r$ and~$\rho_j(\zeta_q)=\zeta_q^j$.  The map $(\ZZ/q\ZZ)^*\longrightarrow G$ given by $j\mapsto \rho_j$ is an isomorphism of groups.
We write the actions of the group rings $\ZZ[\Delta]$ and $\ZZ[G]$ on the multiplicative group
$\QQ(\zeta_r^{},\zeta_q^{})^*$ using exponential notation.

One easily checks the following relations.
$$
\tau(\chi)^{\sigma_i}=\tau(\chi^i),\qquad\hbox{for $i\in(\ZZ/r\ZZ)^*$.}
$$
and
$$
\tau(\chi)^{\rho_j} =\chi(j)^{-1}\tau(\chi),\qquad\hbox{for $j\in(\ZZ/q\ZZ)^*$.}
$$
We write $\overline{\tau(\chi)}$ for the complex conjugate of ${\tau(\chi)}$. For $\chi\not=1$ 
one has 
$$
\tau(\chi)\overline{\tau(\chi)}=q,
$$
showing that $\tau(\chi)$ is an algebraic integer that is only divisible by primes that lie over~$q$.

For our purposes the key property of the Gaussian sums is the following.

\proclaim Proposition 3.1. Let  $q$ be a prime, let $r$ be a positive integer prime to~$q$.
Let $\chi:(\ZZ/q\ZZ)^*\longrightarrow \mu_r$ be
a character  and let $\tau(\chi)$ be the corresponding 
Gaussian sum. Then, for every prime number~$p$ not dividing $qr$ we have
$$
\tau(\chi)^{\sigma_p-p}=\chi^p(p),\quad
\hbox{in the ring $\ZZ[\zeta_q,\zeta_r]/(p)$.}
$$

\smallskip\noindent {\bf Proof.} We have that $\tau(\chi)^p\equiv-\sum_{x\in(\ZZ/q\ZZ)^*}\zeta^{px}_q\chi^p(x)$
modulo the ideal $p\ZZ[\zeta_q^{},\zeta_r^{}]$.
Multiplying by $\chi^p(p)$ and replacing the variable $x$ by~$p^{-1}x$, we get that 
$$\chi^{p}(p)\tau(\chi)^p\equiv
-\chi^{p}(p)\sum_{x\in(\ZZ/q\ZZ)^*}\zeta^{x}_q\chi^p(p^{-1}x)=\tau(\chi^p)\equiv\tau(\chi)^{\sigma_p}\zmod{p}
$$ 
as required.

\medskip 
The cyclotomic primality test proceeds by checking 
the congruence of Proposition~3.1 for suitable  characters~$\chi:(\ZZ/q\ZZ)^*\longrightarrow \mu_r$. 
The next theorem is the 
key ingredient for the cyclotomic primality test.

\proclaim Theorem 3.2. Let $n$ be a natural number. 
Let  $q$ be a prime not dividing~$n$, let $r$ be a power of a prime number~$l$
not dividing~$n$
and let $\chi:(\ZZ/q\ZZ)^*\longrightarrow \mu_r$ be a character.
If
\item{--} for every prime $p$ dividing~$n$ there exists   $\lambda_p$ in the ring $\ZZ_l$ of $l$-adic integers
such that
$$
p^{l-1}= n^{(l-1)\lambda_p},\qquad\hbox{in $\ZZ_l^*$;}
$$
\item{--} the Gaussian sum $\tau(\chi)$  satisfies
$$
\tau(\chi)^{\sigma_n-n}\in\langle\zeta_{r}\rangle,\quad\hbox{in the ring $\ZZ[\zeta_q,\zeta_r]/(n)$},
$$\vskip 0in\noindent
then we have 
$$
\chi(p)=\chi(n)^{\lambda_p}
$$
for every prime divisor $p$ of~$n$.

\smallskip\noindent  Note that  $\lambda_p\in\ZZ_l$ in the first condition is well defined because 
both $n^{l-1}$ and $p^{l-1}$ are congruent to $1\zmod l$.  In addition, $\lambda_p$ is unique. 
When $l$ is odd,  the first condition is equivalent to the condition that
the fraction $(p^{l-1}-1)/(n^{l-1}-1)$ is $l$-integral. 
 In the second condition, we denote by $\langle \zeta_r\rangle$ the cyclic subgroup
of $(\ZZ[\zeta_r^{}]/(n))^*$ of order~$r$ generated by~$\zeta_r$. Note that
the group $\langle \zeta_r\rangle$ is not necessarily equal to
the group  of $r$-th roots of unity in the ring~$\ZZ[\zeta_r^{}]/(n)$.

\medskip\noindent 
{\bf Proof of the theorem.}  We may assume that $\chi$ is a non-trivial character.
By the second condition we have that
$$
\tau(\chi)^{\sigma_n^{-1}n}=\eta\tau(\chi),\qquad\hbox{for some $\eta\in\langle\zeta_{r}\rangle\,\,\subset
\,\,\ZZ[\zeta_q,\zeta_r]/(n)$}.
$$
Note that the operator $\sigma_n^{-1}n\in\ZZ[\Delta]$ has the property that~$\eta^{\sigma_n^{-1}n}=1$.
Therefore, for any integer $L\ge0$, applying  it $(l-1)L$ times leads to the relation
$$
\tau(\chi)^{(\sigma_n^{-1}n)^{(l-1)L}}=\eta^{(l-1)L}\tau(\chi),
\quad\hbox{in the ring $\ZZ[\zeta_q,\zeta_r]/(n)$}.
$$
On the other hand, Proposition~3.1 implies that for any prime divisor $p$ of $n$ we have
$\tau(\chi)^{\sigma_p^{-1}p}=\chi(p)^{-1}\tau(\chi)$ and hence
$$
\tau(\chi)^{(\sigma_p^{-1}p)^{l-1}}=\chi(p)^{1-l}\tau(\chi)\qquad\hbox{in the ring $\ZZ[\zeta_q,\zeta_r]/(p)$.}
$$
Let $l^M$ be the order of the $l$-part of the  finite multiplicative group $(\ZZ[\zeta_q,\zeta_r]/(n))^*$
and let $A$ denote the group $(\ZZ[\zeta_q,\zeta_r]/(n))^*$ modulo $l^M$-th powers.
Let  $L$  be an integer between $0$ and $l^M$ for which $L\equiv \lambda_p\zmod{l^M}$.
Then 
we have $p^{l-1}\equiv n^{(l-1)L}\equiv\zmod{l^M}$ and hence
$(\sigma_n^{-1}n)^{(l-1)L}=\sigma_p^{-1}p$ in the
ring~$(\ZZ/l^M\ZZ)[\Delta]$.
It follows that the left hand sides of the two formulas above are equal in the group~$A$.
Then the same is true for the right hand sides. Since $\tau(\chi)$ is invertible modulo~$p$,
this means
$$
\eta^{(l-1)L}=\chi(p)^{1-l},\qquad\hbox{in the group $A$.}
$$
Since $l-1$ is coprime to the order of $\mu_r$ and since the natural map 
$\langle\zeta_r^{}\rangle\hookrightarrow A$ is injective, this implies
$$
\chi(p)^{-1}=\eta^L=\eta^{\lambda_p},
$$
in the group $\langle\zeta_r\rangle\subset (\ZZ[\zeta_q,\zeta_r]/(n))^*$.
When we multiply the formulas of the first condition for the various prime divisors $p$ of~$n$ 
together, we see that for every positive divisor $d$ of $n$ there exists $\lambda_d\in\ZZ_l$
for which $d^{l-1}= n^{(l-1)\lambda_d}$ in~$\ZZ_l$. We have, of course, $\lambda_n=1$.
From the relation
$\lambda_{dd'}=\lambda_d+ \lambda_{d'}$,
we deduce that $\eta^{\lambda_d}=\chi(d)^{-1}$ for every divisor~$d$ of~$n$.
In particular, we have $\eta=\eta^{\lambda_n}=\chi(n)^{-1}$ and hence 
$$
\chi(p)=\chi(n)^{\lambda_p},
$$
for every prime divisor $p$ of~$n$, as required.
\bigskip\noindent
{\bf Algorithm.} 
The following algorithm is based on Theorem~3.2. 
Suppose we want to prove that a natural number $n$ is prime. 
First determine an integer $R>0$ that has the property that
$$
s\,\,=\,\,\prod_{{q-1|R}\atop{q\,\rm prime}}q
$$ 
exceeds $\sqrt{n}$. At the end of this section we recall that there is a 
constant $c>0$ so that for every natural number $n>16$ there exists an integer $R<(\log\, n)^{c\,\log\,\log\,\log n}$ that has this property. Taking $R$ equal to the product of the first few small prime powers is a good choice.
For all primes $q$ dividing~$s$ and 
for each prime power $r$ that divides $q-1$ exactly, we make sure that ${\rm gcd}(n,qr)=1$ and then check 
the two conditions of Theorem~3.2 for one character of conductor $q$ and order~$r$.
When $n$ passes all these tests, we check for $k=1,\ldots,R-1$ whether the smallest positive residue of~$n^k$ modulo~$s$ divides~$n$.
If that never happens, then $n$ is prime.
\medskip\noindent{\bf Proof of correctness.}
We first note that when $n$ is prime, Proposition~3.1 implies that
it passes all tests.  Conversely, 
suppose that $p\le\sqrt{n}$ is a prime divisor of~$n$.
For every prime $l$ dividing $R$, let $\lambda_p$ be the $l$-adic number that occurs in the first condition of Theorem~3.2.
Let $L\in\{0,1,\ldots,R-1\}$ be the unique integer for which we have
$$L\equiv\lambda_p\zmod{r},
$$
for the power $r$ of $l$ that  exactly divides~$R$.
Theorem~3.2 implies therefore that $\chi(p)=\chi(n)^{L}$ for the  set of characters
of conductor $q$ and order~$r$ for which the conditions of Theorem~3.2 have been checked.  Since
we have $s=\prod_{q-1|R}q$,  the exponent of the group $(\ZZ/s\ZZ)^*$ divides~$R$. 
 Therefore  our set of characters {\it generates} the group of {\it all} characters of~$(\ZZ/s\ZZ)^*$.
It follows that 
$$p\equiv n^{L}\zmod{s}.
$$
Since we have $0<p\le\sqrt{n}<s$, this means that $p$ must actually be {\it equal} to 
the smallest positive residue of~$n^k$ modulo~$s$ for some $k=0,1,\ldots,R-1$.
Since we checked that neither of these numbers divide~$n$, we obtain a contradiction.
It follows that $p$ cannot exist, so that $n$ is necessarily prime.

\medskip

In practice,  checking the first condition of Theorem~3.2 is easy.
When $l\not=2$,  the number $\lambda_p\in\ZZ_l$ of the first condiction
exists if and only if for any prime divisor $p$ of~$n$,
the rational number  $(p^{l-1}-1)/(n^{l-1}-1)$ is $l$-integral. 
Since we have $p^{l-1}\equiv 1\zmod{l}$, this is automatic when we have~$n^{l-1}\not\equiv 1\zmod{l^2}$. 
Given $n$, this  usually holds true  for various prime numbers~$l$. 
Another useful criterion is the following. It can be checked `for free' when one checks
 the second condition of Theorem~3.2.

\proclaim Proposition 3.3. Let $n>1$ be an integer 
and let  $l$ be a prime number  not dividing~$n$.
Then there exists for every prime divisor $p$ of $n$ an exponent $\lambda_p\in\ZZ_l$
for which  
$$p^{l-1}=n^{(l-1)\lambda_p}\quad\hbox{ in~$\ZZ_l^*$},
$$ 
if there exists a prime $q$ not dividing $n$ for which the following holds.
\smallskip
\item{(i)}  ($l\not=2$) for  some power $r>1$ of~$l$ and some
character~$\chi:(\ZZ/q\ZZ)^*\longrightarrow\mu_r$ of order~$r$ the number
$\tau(\chi)^{\sigma_n-n}$
is  a {\it generator} of the cyclic subgroup $\langle\zeta_r\rangle$ of $(\ZZ[\zeta_q, \zeta_r^{}]/(n))^*$.
\item{(ii)} ($l=2$ and $n\equiv1\zmod 4$) we have $\tau(\chi)^{\sigma_n-n}=-1$ for 
the  quadratic character $\chi$ modulo~$q$.
\item{(iii)} ($l=2$ and  $n\equiv3\zmod 4$) and  for  some  character
 $\chi:(\ZZ/q\ZZ)^*\longrightarrow\mu_r$ of 2-power order  $r\ge 4$, the number
 $\tau(\chi)^{\sigma_n-n}$  is  a {\it generator} of the cyclic subgroup 
 $\langle\zeta_r\rangle$ of $(\ZZ[\zeta_q, \zeta_r^{}]/(n))^*$. Moreover,  the Gaussian sum associated to the quadratic character
$\chi^{r/2}$ satisfies  $\tau(\chi^{r/2})^{\sigma_n-n}=-1$ in the ring~$\ZZ[\zeta_q]/(n)$.
 
 \smallskip\noindent{\bf Proof.} Let $p$ be a prime divisor of~$n$ and let $r$ be a power of~$l$.
  As in the proof of Theorem~3.2,
  let $l^M$ denote the order of the $l$-part of the unit group $(\ZZ[\zeta_q,\zeta_{r}]/(p))^*$ and let $A$  be the group
 $(\ZZ[\zeta_q,\zeta_{r}]/(p))^*$ modulo $l^M$-th powers. The latter is a module over the $l$-adic group ring $\ZZ_l[\Delta]$.
 The multiplicative subgroup  
 $\{\sigma_m^{-1}m\in\ZZ_l[\Delta]:m\in\ZZ_l^*\}$ is naturally isomorphic to~$\ZZ^*_l$.
 Therefore, when $l\not=2$, its subgroup $G$ of $(l-1)$-th powers is isomorphic to the additive group~$\ZZ_l$. 
 When $l=2$, this is not true,
 but in that case  the subgroup $G^2$ of squares is isomorphic to~$\ZZ_2$.
  By Proposition~3.1 for any prime $q$ and character $\chi:(\ZZ/q\ZZ)^*\longrightarrow\mu_r$ of order~$r$  we have
  $$\tau(\chi)^{\sigma_p^{-1}p}=\chi(p)^{-1}\tau(\chi),\qquad\hbox{in the group~$A$.}
$$
If  $\tau(\chi)^{\sigma_n-n}$ is  a generator of the group 
$\langle\zeta_{r}\rangle\subset(\ZZ[\zeta_q,\zeta_r^{}]/(n))^*$, then we have
 $$
 \tau(\chi)^{\sigma_n^{-1}n}=\eta\tau(\chi),\qquad\hbox{in the group~$A$.}
 $$
 for some primitive $r$-th root of unity $\eta\in\langle\zeta_{r}\rangle\subset(\ZZ[\zeta_q,\zeta_r^{}]/(n))^*$.
 
 Now we prove {\sl (i)}.  Since $\eta$ is a primitive root, the operator $(\sigma_n^{-1}n)^{l-1}\in\ZZ_l[\Delta]$
 cannot be a `proper' $l$-adic  power of $(\sigma_p^{-1}p)^{l-1}$ in the sense that there cannot 
 exist $\mu\in l\ZZ_l$ for which $(\sigma_n^{-1}n)^{l-1}=(\sigma_p^{-1}p)^{\mu(l-1)}$.
 Since both operators are contained in the pro-cyclic group $G\cong\ZZ_l$, the converse
 must therefore be true: we have $(\sigma_p^{-1}p)^{l-1}=(\sigma_n^{-1}n)^{(l-1)\lambda_p}$ 
 and hence $p^{l-1}=n^{(l-1)\lambda_p}$ for some $\lambda_p\in\ZZ_l$.
 
 To prove~{\sl (ii)}, we observe that the values of $\chi$ are either $1$ or~$-1$.
 Therefore we have $\tau(\chi)^{\sigma_n}=\tau(\chi)$. Since we have
 $\tau(\chi)^2=\chi(-1)\tau(\chi)\overline{\tau(\chi)}=\chi(-1)q$, the condition 
 $\tau(\chi)^{\sigma_n-n}=-1$ means precisely that 
 $$
 (\chi(-1)q)^{(n-1)/2}\equiv -1\zmod{n}.
 $$
 This shows that 
 the 2-parts of the order of $\chi(-1)q\zmod p$ and of $n-1$ are equal. This means that
 $n-1$ divides $p-1$ in the ring of 2-adic integers~$\ZZ_2$. Since $n\equiv 1\zmod 4$, this is equivalent to
 the statement that $p=n^{\lambda_p}$ for some $\lambda_p\in\ZZ_2$.

 To prove~{\sl (iii)}, we note that for $l=2$, the group $G$ that we considered above is not 
 isomorphic to~$\ZZ_2$, but
 the subgroup $G^2$ is. Therefore the arguments of the proof of part~{\sl (i)} only show that
$p^2=n^{2\lambda_p}$  and hence $p=\pm n^{\lambda_p}$ for some  $\lambda_p\in\ZZ_2$.
We show that we have the plus sign. From the relation  $p^2=n^{2\lambda_p}$ we deduce that
$\chi^{-1}(p)^2=\eta^{2\lambda_p}$. Raising this relation to the power $-r/4$, we find
$$
\left({p\over q}\right)\,\,=\,\,\chi^{r/2}(p)\,\,=\,\,\eta^{-r\lambda_p/2}\,\,=\,\,(-1)^{\lambda_p}.
$$
Here we used the usual Legendre symbol to denote the quadratic character $\chi^{r/2}$.
Since $q\equiv 1\zmod 4$, we have $\chi(-1)=1$. Therefore
the second condition  $\tau(\chi^{r/2})^{\sigma_n-n}\equiv-1\zmod{n}$  
says precisely that we have  $q^{(n-1)/2}\equiv -1\zmod n$.
Since $(n-1)/2$ is odd, it follows that
$$
\left({q\over p}\right)\,\,=\,\,\left({{q^{(n-1)/2}}\over p}\right)\,\,=\,\,\left({{-1}\over p}\right).
$$
Since $\chi$ has order at least~4, we have $q\equiv 1\zmod 4$ and hence, by quadratic reciprocity, 
$\left({p\over q}\right)=\left({q\over p}\right)$. The two formulas above imply that
$\left({{-1}\over p}\right)=(-1)^{\lambda_p}$. This means precisely that $p\equiv n^{\lambda_p}\zmod 4$, so that we must have the plus sign, as required.

\medskip
If the number $n$ that is being tested for primality is actually prime, then in each instance the conditions
of Proposition~3.3 are satisfied for a prime $q$ that has the property that
$n$ is not an $l$-th power modulo~$q$. Given $n$, one
encounters in practice for every prime $l$ very quickly such a prime $q$, 
so that the first condition of Theorem~3.2 can be verified.
In the unlikely event that  for some prime~$l$ none of the primes $q$ has this property, 
one simply tests the second condition of Theorem~3.2
for some more primes~$q\equiv 1\zmod l$.

Testing the second condition of Theorem~3.2 is a straightforward computation in the finite ring~$\ZZ[\zeta_q, \zeta_r^{}]/(n)$.
In practice it is important to reduce this to a computation in the much smaller subring~$\ZZ[\zeta_r^{}]/(n)$. This is done by 
using Jacobi sums. 

\medskip
\noindent{\bf Definition.} {\sl Let $q$ be a prime and let 
$\chi,\chi':(\ZZ/q\ZZ)^*\longrightarrow \mu_{r}$ be
two characters.
Then we define the {\it Jacobi sum} $j(\chi,\chi')$ by
$$
j(\chi,\chi')=-\sum_{x\in\ZZ/q\ZZ}\chi(x)\chi'(1-x).
$$
Here we extend $\chi$ and $\chi'$ to $\ZZ/q\ZZ$ by putting~$\chi(0)=\chi'(0)=0$.}

\medskip
The Jacobi sum is an algebraic integer, contained in the  cyclotomic field 
$\QQ(\zeta_r^{})$.
If the characters $\chi,\chi':(\ZZ/q\ZZ)^*\longrightarrow \mu_{r}$  satisfy $\chi\chi'\not=1$,
we have
$$
j(\chi,\chi')={{\tau(\chi)\tau(\chi')}\over{\tau(\chi\chi')}}.
$$
In particular, if $i>0$ is prime to~$r$ and less than the order of~$\chi$, we have
$$
\tau(\chi)^{i-\sigma_i}\,\,=\,\, {{\tau(\chi)^i}\over{\tau(\chi^i)}}\,\,=\,\,\prod_{k=1}^{i-1}j(\chi,\chi^k).
$$
The subgroup of  the $l$-power order roots of unity in~$\overline{\QQ}^*$ is a $\ZZ[\Delta]$-module.
Let $I\subset \ZZ[\Delta]$ be its annihilator. This ideal is generated by the
elements of the form $\sigma_i-i$ with $i\in\ZZ$ coprime to~$l$. Since we have $\tau(\chi)^{\rho_{j}-1}\in\mu_r$ for all $j\not\equiv 0\zmod q$, 
we have
$$
1=\tau(\chi)^{(\rho_{j}-1)x}=\tau(\chi)^{x(\rho_{j}-1)},
\qquad\hbox{for every $x\in I$.}
$$
This shows that $\tau(\chi)^{x}$ and hence that $\tau(\chi)^x$ 
is contained in $\QQ(\zeta_r)$ for every~$x\in\ZZ[\Delta]$. 
This applies in particular to the element $x=\sigma_n-n\in I$. 
It  turns out that it is possible to check the condition of Theorem~3.2
that $\tau(\chi)^{\sigma_n-n}$ is contained in~$\langle\zeta_r\rangle$,
without ever writing down the Gaussian sum $\tau(\chi)\in\ZZ[\zeta_r,\zeta_q]$,
but by doing only computations with Jacobi sums in the ring~$\ZZ[\zeta_r]/(n)$.

When $l$ is odd, the ideal~$I$ generates a {\it principal} ideal 
in the $l$-adic group ring~$\ZZ_l[\Delta]$. It is generated by
any element of the form $\sigma_i-i$ for which $i^{l-1}\not\equiv 1\zmod{l^2}$. We have 
$2^{l-1}\not\equiv1\zmod{l^2}$ for all primes $l<3\cdot 10^9$ except when $l=1093$ or~$3511$.
Therefore we can in practice always use $i=2$. In this case the relevant Jacobi sum is
given by
$$
\tau(\chi)^{\sigma_2-2}\,\,=\,\, {{\tau(\chi)\tau(\chi)}\over{\tau(\chi^2)}}
\,\,=\,\,j(\chi,\chi)\,\,=\,\,-\sum_{x\in\ZZ/q\ZZ}\chi(x(1-x)).
$$
A computation~[\COH, section~9.1.5] shows that we have $\sigma_n-n=\alpha(\sigma_2-2)$ where $\alpha\in\ZZ_l[\Delta]$ 
is given by
$$
\alpha= \sum_{{1\le i<r}\atop{{\rm gcd}(i,r)=1}}\left[{{ni}\over r}\right]\sigma_i^{-1}
$$
times a unit in~$\ZZ_l[\Delta]$. Here $[t]$ denotes the integral part of $t\in{\bf R}$. 
It follows that  in order to verify  that $\tau(\chi)^{\sigma_n-n}$ is contained in 
the group $\langle\zeta_r\rangle$ and to see whether it has order~$r$,  it suffices to evaluate the  product
$$
\prod_{{1\le i<r}\atop{{\rm gcd}(i,r)=1}}j(\chi,\chi)^{\left[{{ni}\over r}\right]\sigma_i^{-1}},
$$
in the ring~$\ZZ[\zeta_r]/(n)$ and check that it is contained in 
the group $\langle\zeta_r\rangle$ and  see whether it has order~$r$. Since the elements in the ring~$\ZZ_l[\Delta]$ map the subgroup $\langle\zeta_r\rangle
\subset(\ZZ[\zeta_r]/(n))^*$ to itself, the fact that we only know the element $\alpha$ up to multiplication by a unit in~$\ZZ_l[\Delta]$ is of no importance.

When $l=2$,  the $\ZZ_l[\Delta]$-ideal generated by~$I$ is {\it not} principal. 
It is generated by the elements $\sigma_3-3$ and $\sigma_{-1}+1$. 
Suppose that the character  $\chi:(\ZZ/q\ZZ)^*\longrightarrow\mu_r$
has 2-power order $r\ge 8$.

When $n\equiv1$ or $3\zmod 8$,  the element
$\sigma_n-n$ is contained in the $\ZZ_l[\Delta]$-ideal 
generated by~$\sigma_3-3$ and we may proceed
as above, replacing the Jacobi sum  by the a product of two Jacobi sums:
$\tau(\chi)^{\sigma_3-3}=j(\chi,\chi)j(\chi,\chi^2)$.
We have $\sigma_n-n=\alpha(\sigma_3-3)$ where $\alpha\in\ZZ_l[\Delta]$ 
is given by  $\alpha= \sum_{i\in E}
\left[{{ni}\over r}\right]\sigma_i^{-1}$ times a unit in~$\ZZ_l[\Delta]$. Here $E$ denotes the set $\{i\in\ZZ:
\hbox{$1\le i<r$ and $i\equiv 1,3\zmod 8$}\}$.
Up to a $\ZZ_l[\Delta]$-automorphism we have
$$
\tau(\chi)^{\sigma_n-n}=\prod_{i\in E}\left(j(\chi,\chi)j(\chi,\chi^2)\right)^{[{{ni}\over r}]
\sigma_i^{-1}},
$$
and this expression involves only elements
in the ring~$\ZZ[\zeta_r]/(n)$.

When $n\equiv 5,7\zmod 8$, we have $\sigma_n-n=-(\sigma_{-n}+n)+(\sigma_{-n}+\sigma_n)$.
Now the element $\sigma_{-n}+n$ is contained in the ideal generated by~$\sigma_3-3$, while
we have $\tau(\chi)^{\sigma_{-n}+\sigma_n}=\tau(\chi^n)\tau(\chi^{-n})=q\chi(-1)$.
In this way one can express  $\tau(\chi)^{\sigma_n-n}$ in a similar way
in terms of elements of the subring~$\ZZ[\zeta_r]/(n)$.
See~[\COH, section~9.1.5] for the formulas

When the order $r$ of the character is 2 or~4,  it is  easier to proceed dircetly.  
When $r=2$, we have $\tau(\chi)^{\sigma_n-n}=(\chi(-1)q)^{(n-1)/2}$ and one should check
that this is equal to $\pm1$ in the ring ${\ZZ/(n)}$. Finally let $r=4$. We have
$\tau(\chi)^{n-\sigma_n}=\left(j(\chi,\chi)^2\chi(-1)q\right)^{(n-1)/4}$ when $n\equiv 1\zmod 4$,
while  $\tau(\chi)^{n-\sigma_n}=j(\chi,\chi)\left(j(\chi,\chi)^2\chi(-1)q\right)^{(n-3)/4}$ when $n\equiv 3\zmod 4$.
In either case, in order to verify

 the second condition of Theorem~2.3, one should check that this number is  a power of $i$ in the ring $\ZZ[i]/(n)$.

\smallskip\noindent{\bf Running time analysis.}
All computations take place in finite rings of the form $\ZZ[\zeta_r]/(n)$, where $r$ divides $R$.
The various summations range over
the congruence classes modulo~$r$ or~$q$.  Both $q$ and $r$ are less than~$R$.
The number of pairs $(q,r)$ involved in the computations
is also at most~$O(R)$.
It follows that  the number of  elementary operations needed to perform the calculations is proportional to $R$ times 
a power of $\log\,n$. Therefore it is important that $R$ is small. On the other hand,
 the size of the  $s$ should be  at least~$\sqrt{n}$. 

By a result in analytic number theory~[\CP,~Thm.~4.3.5]  there is a 
constant $c>0$ so that for every natural number $n>16$ there exists an integer $R<(\log\, n)^{c\,\log\,\log\,\log n}$ for which $s=\prod_{q-1|R}q$ exceeds $\sqrt{n}$. 
It follows that the algorithm is almost polynomial time. It runs in time $O((\log\, n)^{c'\,\log\,\log\,\log n})$
for some constant $c'>0$.

For instance, for  $n$ approximately 880 decimal digits, a good choice is 
$R=2^4\cdot3^2\cdot5\cdot7\cdot11\cdot13\cdot17\cdot19$, because then we have $s>10^{441}$.
H.W.~Lenstra proposed a  slight modification of the cyclotomic test,
that  allows one to efficiently test integers satisfying 
$n<s^3$ rather than $n<s^2$, for primality. See~[\LEN,~Remark~8.7] and 
[\LRES] for this important practical improvement.

\beginsection 4. The elliptic curve primality test.

The elliptic curve primality test, proposed by A.O.L. Atkin in 1988, is
one of the most powerful primality tests that is used in practice~[\MOR].
In order to explain its principle, we first consider a multiplicative group version of the test.

\proclaim Theorem 4.1. Let $n>1$ be a natural number and suppose that there is an element $a\in\ZZ/n\ZZ$
and an exponent $s>0$ satisfying 
$$\eqalign{a^s&=1;\cr
 a^{s/q}-1&\in(\ZZ/n\ZZ)^*,\qquad\hbox{for every prime divisor $q$ of~$s$.}
\cr}
$$ 
Then any prime dividing $n$ is congruent to $1\zmod s$. In particular, if $s>\sqrt{n}$, then $n$ is prime.

\smallskip\noindent{\bf Proof.}  Let $p$ be a prime divisor of~$n$. Then the image of $a$ in $\ZZ/p\ZZ$ is a unit of order~$s$. Indeed, $a^s\equiv1\zmod p$  while $a^{s/q}\not\equiv 1\zmod{p}$  for every prime divisor $q$ of~$s$. Therefore $s$ divides the order  of $(\ZZ/p\ZZ)^*$. In other words, $p\equiv 1\zmod{s}$, as required.
Since a composite $n$ has a prime divisor $p\le n$, the second statement of the theorem is also clear.
Therefore the theorem follows.

\medskip
In applications, $s$ is a divisor of~$n-1$ and the element $a\in\ZZ/n\ZZ$ is the $(n-1)/s$-th power of a  randomly selected element. In order to test the condition that   $a^{s/q}-1\in(\ZZ/n\ZZ)^*$ for every prime divisor $q$ of~$s$,
one evaluates the powers $b=a^{s/q}$ in the ring $\ZZ/n\ZZ$ and then checks that~${\rm gcd}(n,b-1)=1$.
In order to do this, one needs  to know all prime divisors $q$  of~$s$. On the other hand, 
$s$ needs to be large!. Indeed, in order to conclude that 
$n$ is prime, one needs that $s>\sqrt{n}$. In practice, $s$ a completely factored divisor of~$n-1$. If $n$ is large, computing such a divisor of $n-1$ is usually very time consuming.
Therefore, only rarely a large number $n$  is proved prime by a direct application of this theorem.

Occasionally however, it may happen that one can compute a  divisor $r>1$ of $n-1$ that has the property that $s=(n-1)/r$ is {\it probably} prime. In practice, $r$ is the product of the small prime divisors of $n-1$ that one is able to find in a reasonable short time. Therefore $r$ is rather small. Its cofactor $s$  is much larger.
If, by a stroke of luck, the number $s$  happens to pass some probabilistic primality test and one is confident that $s$ is prime, then Theorem~4.1 reduces the problem of proving the primality of~$n$ to proving 
the primality of~$s$, which is at most half the size of~$n$
and usually quite a bit smaller.  Indeed, pick a random $x\in\ZZ/n\ZZ$ and compute $a=x^r$. 
With very high probability we have $a^s\equiv 1\zmod n$ and~$a-1\in(\ZZ/n\ZZ)^*$.
Since $s>\sqrt{n}$, Theorem~4.1 implies  that $n$ is prime {\it provided that} 
the smaller number $s$ is prime. 
However, the chance that $n-1$ factors this way is on the average~$O({1\over{\log\,n}})$.  
Therefore any attempt to proceed in some kind of inductive way, has only a very slight chance of succeeding.

Elliptic curves provide a way out of this situation. The main point is that for prime $n$
there are {\it many} elliptic curves $E$ over~$\ZZ/n\ZZ$ and 
 the orders  of the groups $E(\ZZ/n\ZZ)$ are rather uniformly distributed in the interval $(n+1-2\sqrt{n},n+1+2\sqrt{n})$.
In 1986, S.~Goldwasser and J.~Kilian~[\GK] proposed a primality test based on the principle of Theorem~4.1 
and on a deterministic  polynomial time algorithm to determine 
the number of points on an elliptic curve  over a finite field~[\SCH].
The  running time of their probabilistic algorithm is polynomial time if one assumes a 
certain unproved assumption on the distribution of prime numbers
in short intervals.  Some years later, L.~Adleman and  M.-D.~Huang eliminated the assumption,
by proposing a probabilistic test~[\AH]  involving abelian varieties of dimension~2.  Both tests
are  of theoretical rather than practical value.  By now, even from a theoretical point of view 
they have been superseded by the much simpler  polynomial time deterministic algorithm
explained in section~2.

The key result is the following elliptic  analogue of Theorem~4.1.

\proclaim Theorem 4.2. Let $n>1$ be a natural number and let $E$ be an elliptic curve over $\ZZ/n\ZZ$.
Suppose that there is a point $P\in E(\ZZ/n\ZZ)$ and an integer $s>0$ for which 
$$\eqalign{sP&=0,\quad\hbox{in $E(\ZZ/n\ZZ)$};\cr
\hbox{${s\over q}$}P&\not= 0,\quad\hbox{in $E(\ZZ/p\ZZ)$ for any prime divisor $p$ of~$n$.}\cr}
$$
Then every prime $p$ dividing~$n$ satisfies $\#E(\ZZ/p\ZZ)\equiv 0\zmod{s}$. In particular,  if $s>({\root 4\of n}+1)^2$, then $n$ is prime.

\smallskip\noindent{\bf Proof.}  Let $p$ be a prime divisor of~$n$. Then the image of the point $P$ in $E(\ZZ/p\ZZ)$ has order~$s$. This implies that $\#E(\ZZ/p\ZZ)\equiv 0\zmod{s}$. By Hasse's Theorem, we have that 
$\#E(\ZZ/p\ZZ)\le(\sqrt{p}+1)^2$. Therefore, if $s>({\root 4\of n}+1)^2$, we have that 
$$(\sqrt{p}+1)^2\quad\ge\quad\#E(\ZZ/p\ZZ)\quad\ge \quad s\quad\ge\quad({\root 4\of n}+1)^2
$$ 
and hence $p>\sqrt{n}$. If $n$ were composite, it would have a prime divisor $p\le \sqrt{n}$. 
We conclude that  $n$ is prime as required.

\medskip
The algorithm  reduces the problem of proving the primality of~$n$, to the problem 
of proving that a  smaller number is prime as follows. Given a probable prime number $n$, one randomly 
selects elliptic curves $E$ over $\ZZ/n\ZZ$ and determines the order of the group $E(\ZZ/n\ZZ)$
until one finds a curve for which $\#E(\ZZ/n\ZZ)$ is of the form $r\cdot s$, 
where $s$ 
is a {\it probable} prime number satisfying $s>({\root 4\of n}+1)^2$.
In order to apply Theorem~4.2, one selects a random point $Q\in E(\ZZ/n\ZZ)$ and  computes $P=rQ$.
One checks that $sP=0$ in $E(\ZZ/n\ZZ)$ and that $P\not=0$ in $E(\ZZ/p\ZZ)$
for every prime dividing~$n$. If one works with projective coordinates satisfying a Weierstrass equation,
then the latter  simply means that the gcd of $n$ and the $z$-coordinate of $P$  is equal to~1.
Theorem~4.2 implies then  that $n$ is prime {\it if} $s$ is prime.

In practice, one computes $\#E(\ZZ/n\ZZ)$ {\it under the assumption that $n$ is prime}. Then one
attempts to factor the order of the group $E(\ZZ/n\ZZ)$ by means of a simple trial divison
algorithm or another method that finds small prime factors quicker than larger ones, like Lenstra's Elliptic Curve Method~[\ECM]. Let $r$ be the product of these small prime factors. When $\#E(\ZZ/n\ZZ)$ factors
as  a product $r\cdot s$ with $s$ a probable prime, 
it is in practice not a problem to verify the conditions of Theorem~4.2 for some
randomly selected a point $P$. That's because $n$ is probably prime. But
we do not need to know this in order to apply Theorem~4.2.

Just as in the multiplicative case discussed above,  this computation usually does not work out
when $n$ is large. Typically
one only succeeds in computing 
a small completely factored factor $r$ of $\#E(\ZZ/n\ZZ)$ whose cofactor $s$ is {\it not} prime, but 
cannot be factored easily. In that case one discards the curve $E$, randomly selects another one and tries again. 
Since the curves $E$ are rather uniformly distributed with respect to the number of points
in~$\#E(\ZZ/n\ZZ)$,  the number of attempts one needs to make  before one 
encounters a {\it prime} cofactor~$s$, is expected to be $O(\log\,n)$.
In the unlikely event that one is able to factor $\#E(\ZZ/n\ZZ)$ completely or that 
one has $s<({\root 4\of n}+1)^2$, one is also satisfied. If this happens, 
one can switch the roles of $r$ and $s$ and almost certainly apply Theorem~4.2.

Atkin turns the test of Goldwasser and Kilian into a {\it practical} test by 
selecting the elliptic curves $E$ in the algorithm above more carefully~[\AM]. Atkin
considers suitable
elliptic curves over the complex numbers with {\it complex multiplication} (CM) by 
imaginary quadratic orders of relatively small discriminant. He reduces the curves modulo~$n$ and uses
only these in his primality proof.
The main point is that it is not only theoretically, but also {\it in practice} 
very  easy to count the number of points on these elliptic curves modulo~$n$.
The resulting test is in practice very efficient, but its running time
is very difficult to analyze rigorously, even assuming various conjectures 
on the distributions of smooth numbers and prime numbers. We merely outline the algorithm.

Given $n$, Atkin  first searches for imaginary quadratic  integers 
$\varphi$, for which the following two conditions hold.
$$\eqalign{N(\varphi)&=n,\cr
N(\varphi-1)&=r\cdot s.\cr}
$$
where  we have $r>1$ and   where  $s$ 
satisfies $s>({\root 4\of n}+1)^2$ and is probably prime,
in the sense that it passes a probabilistic primality test.
Here $N(\alpha)$ denote the {\it norm} of an imaginary quadratic number~$\alpha$.

The theory of complex multiplication guarantees the existence of an elliptic curve $E$ over~$\bf C$ 
with endomorphism ring isomorphic to the ring of integers of the
imaginary quadratic field~$\QQ(\varphi)$. Moreover, if $n$ is prime, the characteristic polynomial
of the Frobenius endomorphism of the reduced curve $E\zmod n$  is equal to the minimum polynomial of~$\varphi$.
The number of points in $E(\ZZ/n\ZZ)$ is equal to $N(\varphi-1)=r\cdot s$. Therefore one may
apply Theorem~4.2 to some randomly selected point and conclude that $n$ is prime when $s$ is.
We first explain how to compute suitable imaginary quadratic  integers 
$\varphi$ and then how to compute the corresponding elliptic curves.

If $n$ is prime,  an imaginary quadratic field $F$ contains an element $\varphi$ with $N(\varphi)=n$ if and only if $n$ factors
as a product of two {\it principal} prime ideals in the ring of integers $O_F$ of~$F=\QQ(\varphi)$. 
The probability  that this happens is $1/2h$ where $h$ is the class number of~$O_F$. 
Therefore in practice one first considers all imaginary quadratic fields with class number~$h=1$, then the ones with class number~$h=2$, $\ldots$, etc.  First one checks whether or not $n$ splits in~$F$. If $n$ is prime, this happens if and only if the discriminant $\Delta_F$ is a square modulo~$n$. If $n$ splits, one sees whether it is a product of two prime {\it principal} ideals.
To do this one computes a square root $z$ of $\Delta_F$
modulo~$n$. Then the ideal $I$
generated by $n$ and $z-\sqrt{\Delta_F}$  is a prime divisor of~$n$.
To check that it is principal, one emplys a lattice reduction algorithm and computes
a shortest vector in the rank~2 lattice  generated by $n$ and $z-\sqrt{\Delta_F}$  in~${\bf C}$. If the shortest vector has norm~$n$, then we take it as our integer~$\varphi$ and we know that $I=(\varphi)$ is principal. If the norm 
of the shortest vector is not equal to~$n$,
then the ideal $I$ is not principal and there does not exist an algebraic integer $\varphi\in F$ with $N(\varphi)=n$.
In this case we cannot make use of the elliptic curves that have complex multiplication by the ring of integers of~$F$.

In practice one first computes a `chain' of probable prime numbers $n=N(\varphi)$ with $N(\varphi-1)=r\cdot s$ as above, 
with the property that the primality of one number in the chain, implies the primality of the next one.  
The verifications of the condition of Theorem~4.2 for 
the associated elliptic curves $E$ are not expected to pose any problems and are performed
after a suitable chain has been found. 
Computing the chain is a rather unpredictable enterprise, since it depends on how lucky one
is with the attempts to factor the order of the groups~$E(\ZZ/n\ZZ)$. It may involve some backtracking in a tree
of probable primes. We leave this to  the imagination of the reader.

We explain how to compute the
elliptic curves $E$ over $\ZZ/n\ZZ$ from the quadratic integers~$\varphi$.
The $j$-invariants of elliptic curves over~${\bf C}$ that admit
complex multiplication by the ring of integers of $F=\QQ(\varphi)$ are algebraic integers contained 
in the Hilbert class field of~$F$. 
The $j$-invariant of one such curve is given by
$$
j(\tau)={{\left(1+240\sum_{k=1}^{\infty}\sigma_3(k)q^d\right)^3}\over{q\prod_{k=1}^{\infty}(1-q^k)^{24}}},
$$
where we have $q=e^{2\pi i\tau}$ and where $\tau\in{\bf C}$ has positive imaginary part and has the property that the
ring $\ZZ+\ZZ\tau$ is isomorphic to the ring of integers of~$\QQ(\varphi)$.  The function $\sigma_3$ is given by $\sigma_3(m)=\sum_{d|m}d^3$.
The conjugates of $j(\tau)$  conjugates are given by
$j({{\tau +b}\over{a}})$ for suitable integers~$a,b$. One computes approximates to these numbers and then 
computes the  coefficients of
the minimum polynomial of $j(\tau)$. This polynomial is contained in~$\ZZ[X]$ and has huge coefficients.
Therefore one rather works with  modular functions that are contained in 
extensions of moderate degree $d$ (usually $d=12$ or 24) of the function field ${\bf C}(j)$. The coefficients of these modular 
functions are much smaller. Typically their logarithms are $d$ times smaller~[\AM]. 

If $n$ is prime, it  splits by construction completely in the Hilbert class field~$H$. We compute a root of the minimal polyniomial of $j(\tau)$ in~$\ZZ/n\ZZ$ and call it~$j$. From this we compute a Weierstrass equation 
of an elliptic curve $E$ over~$\ZZ/n\ZZ$ with $j$-invariant equal to~$j$.
We perform all necessary  computations as if $n$ were prime.  Since $n$ probably is prime, they will be successful.
If $n$ is prime, then we have $\#E(\ZZ/n\ZZ)=N(\zeta\varphi-1)$ for some root of unity $\zeta\in\QQ(\varphi)$.
If $\zeta\not=1$, we `twist' the curve~$E$ so that we have  
$\#E(\ZZ/n\ZZ)=N(\varphi-1)=r\cdot s$. Usually, we have $\zeta\in\{\pm 1\}$. The exceptions are 
the fields $F=\QQ(i)$ and $F=\QQ(\sqrt{-3})$, in which case there  are 4 and 6 roots of unity respectively.

\bigskip\bigskip\bigskip\bigskip\noindent
\bibliography

\item{[\APR]} Adleman, L., Pomerance C. and Rumely, R.: On distinguishing prime numbers from composite numbers, {\sl Annals of Math.} {\bf 117} (1983) 173--206.
\item{[\AH]}  Adleman, L. and Huang, M.-D.: {\sl Primality Testing And Two Dimensional Abelian Varieties Over Finite Fields}, Lecture Notes In Mathematics {\bf 1512}, Springer Verlag 1992. 
\item{[\AKS]} Agrawal, M., Kayal, N. and Saxena, N.: Primes is in P,  {\sl Annals of Math.}, to appear.
\item{[\AR]} Artjuhov, M.: Certain criteria for the primality of numbers connected with the little 
Fermat theorem (in Russian), {\sl Acta Arithmetica} {\bf 12} (1966/67) 355--364
\item{[\AM]} Atkin, A. and Morain, F.: Elliptic curves and primality proving, {\sl Math. Comp.} {\bf 61} (1993), 29--68.
\item{[\BA]} Bach, E.: Explicit bounds for primality testing and related problems, {\sl Math. Comp.} {\bf 55} (1990) 355--380.
\item{[\COH]} Cohen, H.: {\sl A Course in Computational Algebraic Number Theory},  Graduate 
Texts in Mathematics {\bf 138}, Springer-Verlag, Berlin 1993.
\item{[\CP]} Crandall,~R. and Pomerance, C.: {\sl Prime Numbers;\ a computational perspective}, Springer Verlag, New York 2001.
\item{[\GK]} Goldwasser, S. and Kilian, J.: Almost all primes can be quickly certified. In {\sl Proc. 18th annual ACM Symposium on the theory of computing} (1986), 316--329.
\item{[\KP]} Konyagin, S. and Pomerance, C.: On primes recognizable in polynomial time, in: {\sl
The Mathematics of Paul Erd\H os}, {\bf I}, vol. 13 of {\sl Algorithms and Combinatorics}, 176--198. Springer-Verlag 1997.
\item{[\LANG]} Lang, S.: {\sl Cyclotomic fields}, Graduate Texts in Math.  {\bf 59}, Springer-Verlag, New York 1978.
\item{[\LENMIL]} Lenstra,~H.W.:  Miller's primality test, {\sl Inform. Process. Lett.} {\bf 8} (1979) 86--88. 
\item{[\LEN]} Lenstra,~H.W.: Primality testing algorithms (after Adleman, Rumely and Williams), In S\'em. Bourbaki, Exp.~576, Springer Lecture Notes in Math. {\bf 901}, Springer-Verlag 1981.
\item{[\LRES]} Lenstra,~H.W.: Divisors in residue classes, {\sl Math. Comp.} {\bf 42} (1984) 331--334.
\item{[\ECM]} Lenstra,~H.W.: Factoring integers with elliptic curves, {\sl Annals of Math.} {\bf 126} (1987) 649--673.
\item{[\LP]} Lenstra,~H.W. and Pomerance C.: Primality testing with Gaussian periods, to appear.
\item{[\MIH]} Mih\u ailescu, P.: Cyclotomic primality proving --- Recent developments, {\sl
Proceedings of  ANTS III, Portland, Oregon}, Lecture Notes in Computer Science {\bf 1423} (1998) 95--111.
\item{[\MIL]} Miller, G.: Riemann's hypothesis and tests for primality, {\sl J. Comput. System Sci.} {\bf 13} (1976) 300--317.
\item{[\MOR]} Morain, F.: Primality proving using elliptic curves: An update, {\sl
Proceedings of  ANTS III, Portland, Oregon}, Lecture Notes in Computer Science {\bf 1423} (1998) 111--127.
\item{[\RAB]} Rabin, M.: Probabilistic algorithm for testing primality, {\sl J. Number Theory}  {\bf 12} (1980) 128--138.
\item{[\SCH]} Schoof, R.: Elliptic curves over finite fields and the computation of square roots mod~$p$, {\sl Math. Comp.} {\bf 44} (1985) 483--494.
\item{[\WASH]} Washington, L.: {\sl  Introduction to cyclotomic fields} 2nd edition, Graduate Texts in Math. {\bf 83}, Springer-Verlag, New York 1997.

\bye